\documentclass{amsart}

\usepackage{amsfonts}
\usepackage{graphicx}
\usepackage{stmaryrd}
\usepackage{amssymb}
\usepackage{extarrows}
\usepackage{mathrsfs}
\usepackage{amsthm}
\usepackage{amsmath}
\usepackage{bm}
\usepackage{color}

\newtheorem{theorem}{Theorem}[section]
\newtheorem{corollary}[theorem]{Corollary}
\newtheorem{lemma}[theorem]{Lemma}
\newtheorem{proposition}[theorem]{Proposition}

\newtheorem{remark}{Remark}

\numberwithin{equation}{section}

\newtheorem*{thma}{Theorem A}
\newtheorem*{thmb}{Theorem B}

\def\a{\alpha}

\def\C{\mathbb{C}}
\def\B{\mathbb{B}}

\def\Bn{\mathbb{B}_n}

\def\Cn{\mathbb{C}^n}

\begin{document}

\title[$L^p-L^q$ boundedness of Forelli-Rudin type operators]
{$L^p-L^q$ boundedness of Forelli-Rudin type operators on the unit ball of $\mathbb{C}^n$}

\author{Ruhan Zhao}
\address{Department of Mathematics\\
SUNY Brockport\\
Brockport, NY 14420}
\email{rzhao@brockport.edu}
\thanks{The first author is partially supported by
National Natural Science Foundation of China (No. 11720101003).}

\author{Lifang Zhou}
\address{Department of Mathematics\\
Huzhou University\\
Huzhou, Zhejiang 31300, People's Republic of China}
\email{lfzhou@zjhu.edu.cn}
\thanks{The second author is supported by the National Natural Science
Foundation of China (No. 11801172,  No. 11771139 and No. 12071130).}



\begin{abstract}
We completely characterize $L^p-L^q$ boundedness of two classes
of Forelli-Rudin type operators on the unit ball of $\mathbb{C}^n$
for all $(p, q)\in [1, \infty]\times [1, \infty]$.
The results are not only a complement to some previous results on
Forelli-Rudin type operators
by Kures and Zhu in 2006 and the first author in 2015, but also a high dimension extension
of some results by Cheng, Fang, Wang and Yu in 2017.
\end{abstract}

\maketitle

\section{Introduction}

In this paper we investigate $L^p$-$L^q$ boundedness of two classes of
Forelli-Rudin type operators defined in \cite{Zhu4}.
For any $a, b, c\in \mathbb{R}$, these integral operators are defined by
\[
T_{a, b, c}f(z)=(1-|z|^2)^a\int_{\Bn}\frac{(1-|w|^2)^b}{(1-\langle z, w\rangle)^c}f(w)dv(w),
\]
and
\[
S_{a,b, c}f(z)=(1-|z|^2)^a\int_{\Bn}\frac{(1-|w|^2)^b}{|1-\langle z, w\rangle|^c}f(w)dv(w),
\]
where $\Bn$ denotes the open unit ball in $\Cn$, and $dv(z)$ is
the normalized volume measure on $\Bn$ such that $v(\Bn)=1$.
When $a=0, b=0, c=n+1$, the operator $T_{0, 0, n+1}$ is the holomorphic Bergman projection on $\Bn$,
and when $a=n+1, b=0, c=2(n+1)$, the operator $S_{n+1, 0, 2(n+1)}$ is the Berezin transform on $\Bn$.
The Bergman projection and the Berezin transform are two of the foundation stones
for the theory of function spaces and the theory of operators, see \cite{ZhuOp}.

Forelli and Rudin initiates the study for boundedness of these operators in \cite{FR},
in which they characterized boundedness of $T_{0,b,n+1+b}$ on $L^p$.
In 1979, Kolaski \cite{KCJ} extended Forelli-Rudin's result from the point of view
of weighted Bergman projection.
In 1991, Zhu \cite{Zhu0} characterized the boundedness of $T_{a, b, c}$ and $S_{a, b, c}$ on $L^p$
with $c=n+1+a+b$ when $n=1$.
In 2006, Kures and Zhu \cite{Zhu4} extended Zhu's result to any value of $c$ except that
$c$ is $0$ or a negative integer.
The authors in this paper respectively observed that the extra requirement for $c$ can be removed
when they studied the norm estimate for $S_{a, b, c}$  in \cite{Zhoulu1}
and the $L^p$-$L^q$ boundedness for $T_{a, b, c}$ and $S_{a, b, c}$ in \cite{Zh}.

For $(p,q)\in [1,\infty]\times [1,\infty]$, the first author \cite{Zh}
in this paper characterized $L^p$-$L^q$ boundedness of $T_{a, b, c}$ and $S_{a,b,c}$
when $1<p\leq q<\infty$ and $1=p\leq q<\infty$.
The range of $p$ and $q$ above forms the red region of the following
$(1/p, 1/q)$-graph in Figure 1.
\begin{center}\setlength{\unitlength}{0.85cm}
\begin{picture}(10,9)
\put(2,2){\includegraphics[width=6cm]{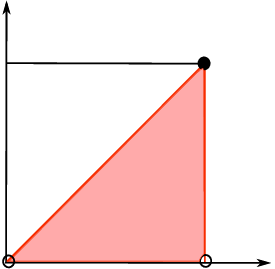}}
\put(2,1.6){$0$}\put(7.1,1.6){$1$}
\put(1.9,7.1){$1$}\put(7.4,7.4){$(1,1)$}
\put(1.8,8.5){$\frac 1q$}\put(8.5,1.6){$\frac 1p$}
\put(1.2,1){$(1/p,  1/q)$-graph (the discussed cases in \cite{Zh})}
\put(4,0.2){Figure 1}
\end{picture}
\end{center}

In \cite{ChFangWY}, Cheng, Fang, Wang and Yu studied the $L^p$-$L^q$ boundedness of the following integral operators
\[
K_cf(z)=\int_{\mathbb{D}} \frac{f(w)}{(1-z\bar{w})^c}\,dA(w)
\]
on the unit disk $\mathbb{D}$ in $\mathbb{C}$, where $dA$ is the normalized area measure of $\mathbb{D}$.
It is clear that $K_c$ is $T_{0, 0, c}$ when $n=1$.
However, in \cite{ChFangWY} they examined the operator $K_c$ from a different angle,
namely, by fixing $c$ and looking at the ranges of $p$ and $q$ such that $K_c:L^p\to L^q$ is bounded.

The main purpose of this paper is to characterize $L^p$-$L^q$ boundedness
of these two classes of Forelli-Rudin type operators in the remaining cases.
As an application,
we will extend the main results, Theorems 1-3 in \cite{ChFangWY},
to the case of weighted Lebesgue spaces on the unit ball.
We note here that our methods are different from the methods in \cite{ChFangWY}.

For $1\leq p<\infty$ and $-1<\a<\infty$, let
$L^{p}_{\a}:=L^{p}(\Bn,\,dv_{\a})$ denote the weighted Lebesgue space
which contains measurable functions $f$ on $\Bn$
such that
$$
\|f\|_{p,\a}=\left(\int_{\Bn}|f(z)|^p\,dv_{\a}(z)\right)^{1/p}<\infty,
$$
where $dv_{\a}(z)=c_{\a}(1-|z|^2)^{\a}\,dv(z)$ and $c_{\a}=\Gamma(n+\alpha+1)/[n!\Gamma(\alpha+1)]$
is the normalized constant such that $v_{\a}(\Bn)=1$.
Let $L^\infty$ be the space of bounded measure function on $\Bn$ such that
\[
\|f\|_\infty=\operatorname*{ess\,sup}_{z\in \Bn}|f(z)|<\infty.
\]
Note that $L_\alpha^\infty=L^\infty$ for any $\a>-1$.

We list the main results in \cite{Zh} as follows. The original results in \cite{Zh} are for operators $T_{0, b, c}$ and $S_{0, b, c}$.
However, for $\alpha>-1$ and $\beta>-1$, it is not difficult to see that $T_{a, b, c}\ (S_{a, b, c})$ is bounded from $L_\alpha^p$ to $L_\beta^q$ if and only if  $T_{0, b, c}\ (S_{0, b, c})$ is bounded from $L_\alpha^p$ to $L_{\beta+aq}^q$ for $1\leq q<\infty$.
The reason we state these results in terms of $T_{a, b, c}$ or $S_{a, b, c}$ is that in this paper we are also considering the case when $q=\infty$, in which the situation is a little subtle.

\begin{thma}\label{thm:pq1}
Suppose $1<p\leq q<\infty$.
Then the following conditions are equivalent.
\begin{itemize}
\item[(i)] The operator $S_{a,b, c}$ is bounded from $L^p_\alpha$ to $L^q_\beta$.
\item[(ii)] The operator $T_{a,b,c}$ is bounded from $L^p_\alpha$ to $L^q_\beta$.
\item[(iii)] The parameters satisfy
\[
\bigg\{
\begin{aligned}
&-qa<\beta+1, \quad (\alpha+1)<p(b+1), \\
&c\leq n+1+a+b+\frac{n+1+\beta}q-\frac{n+1+\alpha}p.
\end{aligned}
\]
\end{itemize}
\end{thma}

\begin{thmb}\label{thm:pq2}
Suppose $1=p\leq q<\infty$.
Then the following conditions are equivalent.
\begin{itemize}
\item[(i)] The operator $S_{a, b, c}$ is bounded from $L^1_\alpha$ to $L^q_\beta$.
\item[(ii)] The operator $T_{a,b, c}$ is bounded from $L^1_\alpha$ to $L^q_\beta$.
\item[(iii)] The parameters satisfy
\[
\bigg\{
\begin{aligned}
&-qa<\beta+1, \quad \alpha<b, \\
&c\leq a+b-\alpha+\frac{n+1+\beta}q,
\end{aligned}
\quad
\mathrm{or}
\quad
\bigg\{
\begin{aligned}
&-qa<\beta+1, \quad \alpha= b, \\
&c<a+\frac{n+1+\beta}q.
\end{aligned}
\]
\end{itemize}
\end{thmb}

Our main results are the following.

\begin{theorem}\label{thm:pq3}
Suppose $1\leq q<p\leq\infty$.
Then the following conditions are equivalent.
\begin{itemize}
\item[(i)] The operator $S_{a, b, c}$ is bounded from $L^p_\alpha$ to $L^q_\beta$.
\item[(ii)]The operator $T_{a, b, c}$ is bounded from $L^p_\alpha$ to $L^q_\beta$
\item[(iii)]The parameters satisfy
\[
\bigg\{
\begin{aligned}
&-qa<\beta+1, \quad (\alpha+1)<p(b+1), \\
&c<n+1+a+b+\frac{1+\beta}q-\frac{1+\alpha}p.
\end{aligned}
\]
When $p=\infty$, these conditions should be interpreted as
\[
\bigg\{
\begin{aligned}
&-qa<\beta+1,\quad b>-1, \\
&c<n+1+a+b+\frac{\beta+1}q.
\end{aligned}
\]
\end{itemize}
\end{theorem}

\begin{theorem}\label{thm:thm4}
Suppose $p=1, q=\infty$. Then the following conditions are equivalent.
\begin{itemize}
\item[(i)] The operator $S_{a, b, c}$ is bounded from $L_\alpha^1$ to $L^\infty$.
\item[(ii)] The operator $T_{a, b, c}$ is bounded from $L_\alpha^1$ to $L^\infty$.
\item[(iii)] The parameters satisfy
\[
\bigg\{
\begin{aligned}
&a\geq 0,\  \alpha\leq b, \\
&c\leq a+b-\alpha.
\end{aligned}
\]
\end{itemize}
\end{theorem}

\begin{theorem}\label{thm:thm5}
Suppose $1< p\leq \infty, \ q=\infty$. Then the following conditions are equivalent.
\begin{itemize}
\item[(i)] The operator $S_{a, b, c}$ is bounded from $L^p_\alpha$ to $L^\infty$.
\item[(ii)] The operator $T_{a, b, c}$ is bounded from $L^p_\alpha$ to $L^\infty$.
\item[(iii)] The parameters satisfy
\[
\bigg\{
\begin{aligned}
&a=0,  \quad \alpha+1<p(b+1),\\
&c<n+1+b-\frac{n+1+\alpha}p,
\end{aligned}
\quad
\mathrm{or}
\quad
\bigg\{
\begin{aligned}
&a>0, \quad (\alpha+1)<p(b+1), \\
&c\leq n+1+a+b-\frac{n+1+\alpha}p.
\end{aligned}
\]
When $p=\infty$, these conditions should be interpreted as
\[
\bigg\{
\begin{aligned}
&a=0, \ b>-1, \\
&c<n+1+b,
\end{aligned}
\quad
\mathrm{or}
\quad
\bigg\{
\begin{aligned}
&a>0, \ b>-1, \\
&c\leq n+1+a+b.
\end{aligned}
\]
\end{itemize}
\end{theorem}

With these three theorems we have completely characterized the $L^p$-$L^q$
boundedness of Forelli-Rudin type operators $T_{a, b, c}$ and $S_{a, b, c}$
for all $(p,q)\in [1,\infty]\times [1,\infty]$.

There are also some researches for the $L^p$-$L^q$ boundedness of Forelli-Rudin type operators on other domains. For example, the $L^p$-$L^q$-boundedness for Forelli-Rudin type operators induced by holomorphic Bergman kernel on the Siegel upper half space was studied in \cite{LSH,SB}.

The paper is organized as follows. In Section 2,  we collect some preliminary results which will be used later. It is easy to see the  $L^p$-$L^q$ boundedness for $S_{a, b, c}$ implies the $L^p$-$L^q$ boundedness of $T_{a, b, c}$,  thus the proof of (i) $\Rightarrow$ (ii) in Theorem \ref{thm:pq3}-Theorem \ref{thm:thm5} are going to be trivial. It suffices for us to prove (iii)$\Rightarrow$ (i) which are sufficiency parts for our main results, and  (ii) $\Rightarrow$ (iii) which are the necessity parts for our main results. Therefore, Section 3 is devoted to giving the proof for the sufficiency,  and Section 4 is devoted to giving the proof for the necessity.
In Section 5, we apply our main results to extend Theorems 1-3 in \cite{ChFangWY} from the unit disk to the unit ball,
and from Lebesgue spaces to weighted Lebesgue spaces.
We also list the results for two special cases, namely, the weighted Bergman projection and the Berezin transform,
as applications of our main results.

\section{Preliminaries}

\subsection{Some integral formulas}

\begin{lemma}\label{lem:formula}
Suppose $z\in \Bn$, $t>-1$, and $c$ is real. The integral
\[
I_{c, t}(z)=\int_{\Bn}\frac{(1-|w|^2)^t}{|1-\langle z, w\rangle|^{c}}dv(w)
\]
has the following asymptotic behavior as $|z|\to 1^{-}$.
\begin{itemize}
\item[(i)] If $n+1+t-c>0$, then $I_{c, t}(z)\approx1$.\\
\item[(ii)] If $n+1+t-c=0$, then $I_{c, t}(z)\approx\log\dfrac{1}{1-|z|^2}$.\\
\item[(iii)] If $n+1+t-c<0$, then $I_{c, t}(z)\approx(1-|z|^2)^{n+1+t-c}$.
\end{itemize}
where the notation $A \approx B$ means that $A/B$ is bounded above and below by two positive constants  that are independent of $z$.
\end{lemma}
\begin{proof}
See \cite[Propostion 1.4.10]{RW}, \cite[Theorem1.12]{ZhuBn} and \cite[Corollary 2.4]{LC}.
\end{proof}

\begin{lemma}\label{lem:formula2}
Suppose $z\in \Bn$, $1\leq p<\infty$. Then
\begin{equation}\label{eq:integral}
\int_{\Bn}\frac{(1-|w|^2)^t}{|1-\langle z, w\rangle|^s}dv(w)\in L_\alpha^p
\end{equation}
if and only if
\[
\bigg\{
\begin{aligned}
&t>-1, \\
&s<n+1+t+\frac{\alpha+1}p.
\end{aligned}
\]
\end{lemma}
\begin{proof}
First, for a fixed $z \in \Bn$, it is clear that in order that the integral is finite, we must have $t>-1$.
Consider the following two cases.

\emph{Case I:  $s<n+1+t+(\alpha+1)/p$. } Choose an $\epsilon\in \big(0, (\alpha+1)/p\big)$ small enough such that $s\leq n+1+t+(\alpha+1)/p-\epsilon$. Thus let
\[
\tau_1=n+1+t+(\alpha+1)/p-\epsilon-s\geq0.
\]
Therefore, according to lemma \ref{lem:formula}, it follows that
\begin{eqnarray*}
\int_{\Bn}\frac{(1-|w|^2)^t}{|1-\langle z, w\rangle|^s}dv(w)
&=&\int_{\Bn}\frac{(1-|w|^2)^t|1-\langle z, w\rangle|^{\tau_1}}{|1-\langle z, w\rangle|^{n+1+t+(\alpha+1)/p-\epsilon}}dv(w)\\
&\leq&2^{\tau_1}\int_{\Bn}\frac{(1-|w|^2)^t}{|1-\langle z, w\rangle|^{n+1+t+(\alpha+1)/p-\epsilon}}dv(w)\\
&\approx& (1-|z|^2)^{\epsilon-(\alpha+1)/p}\in L_\alpha^p.
\end{eqnarray*}

\emph{Case II: $s\geq n+1+t+(\alpha+1)/p$. } Denote  $\tau_2=s-\big(n+1+t+(\alpha+1)/p\big)\geq 0$.  According to lemma \ref{lem:formula} again,  it follows that
\begin{eqnarray*}
\int_{\Bn}\frac{(1-|w|^2)^t}{|1-\langle z, w\rangle|^s}dv(w)
&=&\int_{\Bn}\frac{(1-|w|^2)^t}{|1-\langle z, w\rangle|^{n+1+t+(\alpha+1)/p+\tau_2}}dv(w)\\
&\approx&(1-|z|^2)^{-(\alpha+1)/p-\tau_2}\notin L_\alpha^p.
\end{eqnarray*}
The proof is complete.
\end{proof}

\subsection{Schur's tests}

Let  $p^\prime$ and $q^\prime$ denote the conjugate indices of $p$ and $q$ for $1\leq p, q\leq \infty$ throughout the paper, respectively. We will assume that $p^\prime=\infty$ for $p=1$, and $p^\prime=1$ for $p=\infty$.
The following basic results are referred as Schur's tests.

\begin{lemma}\cite[Theorem 2.2]{SG}\label{lem:1}
Let $\mu$ and $\nu$ be the positive measure on the space $X$ and let $K(x, y)$ be a nonnegative measurable function on $X\times X$. Let $T$ be the integral operator with kernel $K$ defined as follows.
\[
Tf(x)=\int_XK(x, y)f(y)d\mu(y).
\]
Suppose $1<q\leq p<\infty$. If there exist a constant $C>0$ and a positive $\mu$-measurable function $\varphi$ which is finite $\mu$-almost everywhere and satisfies $S^K\varphi\leq C \varphi$ $\mu$-almost everywhere in $X$, then $T$ is bounded from $L_\mu^p$ to $L_\nu^q$. Moreover,
\[
\|T\|_{L^p_\mu\to L^q_\nu}\leq C^{(p-1)/q}\|\varphi\|_{L_\mu^p}^{(p-q)/q},
\]
where
\[
S^K\varphi(y)=\bigg[\int_X K(x, y)\bigg(\int_X K(x,z)\varphi(z)d\mu(z)\bigg)^{q-1}d\nu(x)\bigg]^{p^{\prime}-1}.
\]
\end{lemma}

\begin{lemma}\cite[Propostion 5.4]{Tao1}\label{lem:shurtestpinfty}
Let $\mu$ and $\nu$ be the positive measure on the space $X$ and let  $K: X\times
X\to \mathbb{C}$. Let $T$ be the integral operator with kernel $K$ defined as follows.
\[
Tf(x)=\int_X K(x, y)f(y)d\mu(y).
\]
Let $1\leq p\leq \infty$. If
$\|K(x, \cdot)\|_{L^{p^\prime}_\mu}$ is uniformly bounded, then $T$
is bounded from $L^p_\mu$ to $L_{\nu}^\infty$. Moreover
\[
\|T\|_{L^p_\mu \to L^\infty_\nu}=\operatorname*{ess\,sup}_{x\in X}\|K(x,\cdot)\|_{L^{p^\prime}_\mu}.
\]
\end{lemma}

\begin{lemma}\cite[Problem 5.5]{Tao1}\label{lem:inftyq}
Let $\mu$ and $\nu$ be the positive measure on the space $X$, and $K(x, y)$ be a nonnegative measurable function on $X\times X$. Let $T$ be the integral operator with kernel $K$ defined as follow.
\[
Tf(x)=\int_X K(x, y)f(y)d\mu(y).
\]
Suppose $1 \leq q\leq \infty$.  Then the operator $T$
is bounded from $L^\infty_\mu$ to $L^q_\nu$ if and only if
\[
\int_X K(x, y)d\mu(y)\in L^q_\nu.
\]
Moveover,
\[
\|T\|_{L^\infty_\mu\to L^q_\nu}=\bigg\|\int_X K(x, y)d\mu(y)\bigg\|_{L^q_\nu}
\]
\end{lemma}

\begin{lemma}\cite[Problem 5.5]{Tao1}\label{lem:p1}
Let $\mu$ and $\nu$ be the positive measure on the space $X$ and $K(x, y)$ be a nonnegative measurable function on $X\times X$. Let $1 \leq p< \infty$. Let $T$ be the integral operator with kernel $K$ defined as follows.
\[
Tf(x)=\int_X K(x, y)f(y)d\mu(y).
\]
Then the operator $T$ is bounded from $L^p_\mu$ to $L^1_\nu$ if and only if
\[
\int_X K(x,y)d\nu(x)\in L^{p^\prime}_{\mu}.
\]
Moreover,
\[
\|T\|_{L^p_\mu\to L^1_\nu}=\bigg\|\int_X K(x, y)d\nu(x)\bigg\|_{L^{p^\prime}_{\mu}}.
\]
\end{lemma}

\section{Proof of sufficiency for Boundedness of $S_{a, b,c}$}

\subsection{Proof of (iii)$\Rightarrow$(i)  for  Theorem \ref{thm:pq3}}

The following three propositions are devoted to giving the
the proof of (iii)$\Rightarrow$(i)  for  Theorem \ref{thm:pq3}.
\begin{proposition}\label{lem:11}
Let $1<q<p<\infty$.  If the parameters satisfy
\[
\bigg\{
\begin{aligned}
&-qa<\beta+1, \quad (\alpha+1)<p(b+1),\\
&c<n+1+a+b+\frac{1+\beta}q-\frac{1+\alpha}p
\end{aligned}
\]
then $S_{a, b,c}$ is bounded from $L^p_\alpha$ to $L^q_\beta$.
\end{proposition}

\begin{proof}
The operator $S_{a, b, c}$ is an integral operator with the kernel
\[
K(z,w)=c_\alpha^{-1}\frac{(1-|z|^2)^a(1-|w|^2)^{b-\alpha}}{|1-\langle z, w\rangle|^c}.
\]
Let $c<n+1+a+b+(1+\beta)/q-(1+\alpha)/p$. Then there exists an  $\epsilon>0$, which is small enough such that
\[
a+\frac{\beta+1}q-\frac{p}{p-q}\epsilon>0, \quad (b+1)-\frac{\alpha+1}p-\frac{q(p-1)\epsilon}{p-q}>0,
\]
and
$$
c\leq n+1+a+b+\frac{1+\beta}q-\frac{1+\alpha}p-\epsilon.
$$
Let
\[
\tau=n+1+a+b+\frac{1+\beta}q-\frac{1+\alpha}p-\epsilon-c.
\]
Then $\tau\ge 0$.
Therefore,
\begin{eqnarray*}
K(z, w)
&=&\frac 1{c_\alpha}\frac{(1-|z|^2)^a(1-|w|^2)^{b-\alpha}|1-\langle z, w\rangle|^\tau}{|1-\langle z, w\rangle|^{ n+1+a+b+(1+\beta)/q-(1+\alpha)/p-\epsilon}}\\
&\leq&\frac{2^\tau}{c_\alpha}\frac{(1-|z|^2)^a(1-|w|^2)^{b-\alpha}}{|1-\langle z, w\rangle|^{ n+1+a+b+(1+\beta)/q-(1+\alpha)/p-\epsilon}}.
\end{eqnarray*}
Thus we only need to prove if
\begin{equation}\label{eq:cc}
c=n+1+a+b+\frac{\beta+1}q-\frac{\alpha+1}p-\epsilon,
\end{equation}
then $S_{a,b,c}$ is bounded from $L^p_\alpha$ to $L^q_\beta$.

In the following proof we fix the parameter $c$ to be given by \eqref{eq:cc}.
To apply Lemma~\ref{lem:1}, we need to find a test function $\varphi$ such that
\begin{equation}\label{eq:sabck}
S_{a,b,c}^K\varphi(u)\leq C\varphi(u)
\end{equation}
 almost everywhere in $\Bn$, where
$$
S_{a, b, c}^K\varphi(u)
=\left[ \int_{\Bn}K(z, u)\left(\int_{\Bn}K(z, w)\varphi(w)dv_\alpha(w)\right)^{q-1}dv_{\beta}(z)\right]^{p^{\prime}-1}
$$
and $C$ is a positive constant.
For $w\in \Bn$, let
\[
\varphi(w)=(1-|w|^2)^{q\epsilon/(p-q)-(\alpha+1)/p}.
\]
From our choice of $\epsilon$ we have that
$$
b+\frac{q\epsilon}{p-q}-\frac{\alpha+1}{p}>-1,
$$
and
\[
n+1+\left(b+\frac q{p-q}\epsilon-\frac{\alpha+1}p\right)-c=-\bigg(a+\frac{\beta+1}q-\frac{p}{p-q}\epsilon\bigg)<0.
\]
Thus, by Lemma~\ref{lem:formula} we have
\begin{eqnarray*}
&~&\int_{\Bn}K(z, w)\varphi(w)\,dv_\alpha(w)\\
&~&\qquad= (1-|z|^2)^a\int_{\Bn}\frac{(1-|w|^2)^{b+ q\epsilon/(p-q)-(\alpha+1)/p}}{|1-\langle z, w\rangle|^c}dv(w)\\
&~&\qquad\approx(1-|z|^2)^{-\big((\beta+1)/q-p\epsilon/(p-q)\big)}.
\end{eqnarray*}
Again, from our choice of $\epsilon$ we have that
\begin{eqnarray*}
&~&\beta+a-\left(\frac{\beta+1}q-\frac{p}{p-q}\epsilon\right)(q-1)+1\\
&~&\qquad=a+\frac{\beta+1}q+\frac{p(q-1)}{p-q}\epsilon>0,
\end{eqnarray*}
and
\begin{align*}
&n+1+\left[\beta+a-\left(\frac{\beta+1}q-\frac{p}{p-q}\epsilon\right)(q-1)\right]-c\\
&\quad=-\bigg(1+b-\frac{\alpha+1}p-\frac{q(p-1)}{p-q}\epsilon\bigg)<0.
\end{align*}
By Lemma~\ref{lem:formula} we obtain that
\begin{eqnarray*}
&~&S_{a, b, c}^K\varphi(u)\\
&~&\qquad\approx\left[(1-|u|^2)^{b-\alpha}\int_{\Bn}
\frac{(1-|z|^2)^{\beta+a-\big((\beta+1)/q-p\epsilon/(p-q)\big)(q-1)}}{|1-\langle z,u\rangle|^c}\,dv(z)\right]^{p'-1}\\
&~&\qquad\approx(1-|u|^2)^{-\big((\alpha+1)/p^\prime-q(p-1)\epsilon/(p-q)\big)(p'-1)}\\
&~&\qquad=(1-|u|^2)^{q\epsilon/(p-q)-(\alpha+1)/p}\\
&~&\qquad=\varphi(u).
\end{eqnarray*}
Hence, we have proved ~\eqref{eq:sabck}.
The proof is complete.
\end{proof}

\begin{proposition}\label{prop:13}
Let $1< p< \infty, \ q=1$.  If $\alpha+1<p(b+1)$, then the operator $S_{a, b, c}$ is bounded from $L^p_\alpha$ to $L_\beta^1$ if and only if
\[
\bigg\{
\begin{aligned}
&-a<\beta+1, \\
&c<n+1+a+b+1+\beta-\frac{1+\alpha}p.
\end{aligned}
\]
\end{proposition}

\begin{proof}
It follows from  Lemma \ref{lem:p1} that $S_{a, b, c}$ is bounded from $L^p_\alpha$ to $L_\beta^1$ if and only if
\begin{align*}
(1-|w|^2)^{b-\alpha}\int_{\Bn}\frac{(1-|z|^2)^a}{|1-\langle z, w\rangle|^{c}}dv_\beta(z)\in L_\alpha^{p^\prime},
\end{align*}
which is equivalent to
\begin{equation}\label{eq:kernelestimate}
\int_{\Bn}\frac{(1-|z|^2)^{a+\beta}}{|1-\langle z, w\rangle|^{c}}dv(z)\in L_{p^\prime(b-\alpha)+\alpha}^{p^\prime}.
\end{equation}
From our condition we know that $p^\prime(b-\alpha)+\alpha>-1$.
Thus the result follows from Lemma~\ref{lem:formula2}.
\end{proof}

\begin{proposition}\label{prop:inftyqs}
Let $p=\infty, \ 1\leq q< \infty$.  If $-qa<\beta+1$, then the operator $S_{a, b, c}$ is bounded from $L^\infty$ to $L_\beta^q$ if and only if
\[
\bigg\{
\begin{aligned}
&b>-1, \\
&c<n+1+a+b+\frac{\beta+1}q.
\end{aligned}
\]
\end{proposition}

\begin{proof}
It follows from Lemma \ref{lem:inftyq} that $S_{a, b, c}$ is bounded from $L^\infty$ to $L_\beta^q$ if and only if
\[
(1-|z|^2)^a\int_{\Bn}\frac{(1-|w|^2)^b}{|1-\langle z,w\rangle|^c}dv(w)\in L_\beta^q,
\]
which is equivalent to
\[
\int_{\Bn}\frac{(1-|w|^2)^b}{|1-\langle z,w\rangle|^c}dv(w)\in L_{\beta+qa}^q.
\]
We immediately get the result from Lemma \ref{lem:formula2} since $\beta+qa>-1$.
\end{proof}

\subsection{Proof of (iii)$\Rightarrow $ (i) for Theorem \ref{thm:thm4} and \ref{thm:thm5}}

The sufficiency for boundedness of $S_{a, b, c}$ when $q=\infty$  will be given in the following proposition by using Lemma \ref{lem:shurtestpinfty}.

\begin{proposition}\label{lem:pinfty}
Let $1\leq p\leq \infty, \ q=\infty$. Let the parameters satisfy one of the following conditions,
\begin{itemize}
\item[(i)]
$p=1$, and
\[
\bigg\{
\begin{aligned}
&a\geq 0,\  \alpha\leq b, \\
&c\leq a+b-\alpha.
\end{aligned}
\]
\item[(ii)]
$1<p<\infty$, and
\[
\bigg\{
\begin{aligned}
&a=0,  \quad \alpha+1<p(b+1),\\
&c<n+1+b-\frac{n+1+\alpha}p,
\end{aligned}
\]
or
\[
\bigg\{
\begin{aligned}
&a>0, \quad (\alpha+1)<p(b+1) \\
&c\leq n+1+a+b-\frac{n+1+\alpha}p.
\end{aligned}
\]

\item[(iii)]
 $p=\infty$, and
\[
\bigg\{
\begin{aligned}
&a=0, \ b>-1, \\
&c<n+1+b,
\end{aligned}
\quad
\mathrm{or}
\quad
\bigg\{
\begin{aligned}
&a>0, \ b>-1, \\
&c\leq n+1+a+b.
\end{aligned}
\]
\end{itemize}
Then $S_{a, b,c}$ is bounded from $L^p_\alpha$ to $L^\infty$.
\end{proposition}

\begin{proof}
We first prove (i).
Let $\tau=a+b-\alpha-c$.
Since $c\leq a+b-\alpha$, we know that $\tau\geq0$.
Since $a\geq 0$ and  $\alpha\leq b$, we have that
\begin{eqnarray*}
&~&\left|\frac{(1-|z|^2)^a(1-|w|^2)^{b-\alpha}}{|1-\langle z, w\rangle|^{c}}\right|\\
&~&\qquad\leq \frac{(1-|z|^2)^a(1-|w|^2)^{b-\alpha}|1-\langle z, w\rangle|^\tau}{|1-\langle z, w\rangle|^{a+b-\alpha}}\\
&~&\qquad\leq 2^\tau (1+|z|)^a(1+|w|)^{b-\alpha}\leq 2^{\tau+a+b-\alpha},
\end{eqnarray*}
for all $z\in \Bn$. Thus
$$
\sup_{z\in \Bn}\|K(z, \cdot)\|_{L^{\infty}}
=c_\alpha^{-1}\sup_{z\in \Bn}\bigg\|\frac{(1-|z|^2)^a(1-|w|^2)^{b-\alpha}}{|1-\langle z, w\rangle|^{c}}\bigg\|_{L^\infty}<\infty.
$$
By Lemma~\ref{lem:shurtestpinfty} we know that  $S_{a, b,c}$ is bounded from $L^1_\alpha$ to $L^\infty$.

Now we prove (ii) and (iii). Let $1< p\leq \infty$.
First, since $\alpha+1<p(b+1)$, we have $(b-\alpha)p^\prime+\alpha>-1$.
Next, from the conditions for $c$ in (ii) and (iii) we can easily see that
$$
[n+1+(b-\alpha)p^\prime+\alpha]-cp'>0
$$
for $a=0$, and
$$
[n+1+(b-\alpha)p^\prime+\alpha]-cp'\geq-ap'
$$
for $a>0$.
By Lemma~\ref{lem:formula} we get that
$$
\sup_{z\in \Bn}\|K(z, \cdot)\|_{L_\alpha^{p^\prime}}
=\sup_{z\in \Bn}(1-|z|^2)^{ap^\prime}\int_{\Bn}\frac{(1-|w|^2)^{(b-\alpha)p^\prime+\alpha}}{|1-\langle z, w\rangle|^{cp^\prime}}dv(w)<\infty.
$$
By Lemma~\ref{lem:shurtestpinfty} we know that $S_{a, b,c}$ is bounded from $L^p_\alpha$ to $L^\infty$.
The proof is complete.
\end{proof}

\section{Proof of necessity for Boundedness of $T_{a, b,c}$ }

We need the following lemmas.

\begin{lemma}\label{lem:b1}
Let $1\leq p, q\leq \infty$. If $T_{a, b,c}$ is bounded from   $L_\alpha^p$ to $L_\beta^q$, then $b>-1$.
\end{lemma}
\begin{proof}
Since $T_{a, b,c}$ is bounded from  $L_\alpha^p$ to $L_\beta^q$,
it follows that $T_{a, b, c}1\in L_\beta^q$. Notice
that
\[
(T_{a, b,c}1)(z)=(1-|z|^2)^a\int_{\Bn}\frac{(1-|w|^2)^b}{(1-\langle z, w\rangle)^c}dv(w).
\]
It is clear that in order that the integral is finite for any fixed $z\in\mathbb{B}_n$, we must have $b>-1$.
\end{proof}

\begin{lemma}\label{lem:aqbeta}
Let $1\leq p\leq \infty$ and $1\leq q<\infty$. Suppose that $T_{a, b, c}$ is bounded from $L^p_\alpha$ to $L^q_\beta$, then $-qa<\beta+1$.
\end{lemma}

\begin{proof}
Let $f_N(z)=(1-|z|^2)^N$, where $N$ is a positive real number, and satisfies  $Np+\alpha>-1$ if $1\leq p<\infty$.
It follows from Lemma \ref{lem:formula} that $f_N\in L_\alpha^p$.
The symmetry of $\B_n$ shows that there is a positive constant $C_N$ such that $T_{a,b,c}f_N(z)=C_N(1-|z|^2)^a$.
Since $T_{a,b,c}$ is bounded from $L_\alpha^p$ to $L_\beta^q$, it follows that $(1-|z|^2)^a\in L_\beta^q$,
which implies  $qa+\beta>-1$, or $-qa<\beta+1$. The proof is complete.

\end{proof}

\begin{lemma}\label{lem:pinfty1a}
Let $1\leq p \leq \infty$ and $q=\infty$. Suppose that $T_{a, b, c}$ is bounded from $L_\alpha^p$ to $L^\infty$, then $a\geq 0$.
\end{lemma}
\begin{proof}
Let $f_N(z)=(1-|z|^2)^N$.
A similar argument to Lemma~\ref{lem:aqbeta} shows  $a\geq 0$ from the fact that $T_{a,b,c}f_N(z)=C_N(1-|z|^2)^a\in L^\infty$.
\end{proof}

\begin{lemma}\label{lem:alphap}
Let $1<p<\infty$ and $1\leq q\leq\infty$. Suppose that $T_{a, b, c}$ is bounded from $L^p_\alpha$ to $L^q_\beta$, then $\alpha+1<p(b+1)$.
\end{lemma}
\begin{proof}
Since $T_{a, b, c}$ is bounded from $L^p_\alpha$ to $L^q_\beta$, it follows that its adjoint operator $T_{a,b,c}^\ast$  is bounded from $(L^q_\beta)^\ast$ to $(L^p_\alpha)^\ast$.
A simple computation shows that, when $1\leq q<\infty$, the adjoint operator $T_{a,b,c}^\ast$ is given by
\begin{equation}\label{tstar}
T_{a,b,c}^\ast f(z)=\frac{c_{\beta}}{c_{\alpha}}(1-|z|^2)^{b-\alpha}\int_{\Bn}\frac{(1-|w|^2)^{a+\beta}}{(1-\langle z, w\rangle)^c}f(w)\,dv(w).
\end{equation}
Thus $T_{a,b,c}^\ast$ is bounded from $L^{q^\prime}_\beta$ to $L^{p^\prime}_\alpha$.
When $q=\infty$, the adjoint operator of $T_{a,b,c}$ is given by
\begin{equation}\label{tstar1}
T_{a, b, c}^\ast f(z)=\frac{1}{c_{\alpha}}(1-|z|^2)^{b-\alpha}\int_{\Bn}\frac{(1-|w|^2)^{a}}{(1-\langle z, w\rangle)^c}f(w)\,dv(w)
\end{equation}
which is bounded from $(L^\infty)^\ast$ to $L^{p^\prime}_\alpha$.
Since $L^1$ injects continuously into $(L^\infty)^\ast$, we get that $T_{a, b, c}^\ast$ is bounded from $L^1$ to $L^{p^\prime}_\alpha$.
A similar argument to Lemma~\ref{lem:aqbeta} with boundedness of $T_{a, b, c}^\ast$ shows that
$p'(b-\alpha)+\alpha>-1$, which is the same as $\alpha+1<p(b+1)$.
\end{proof}

\begin{lemma}\label{lem:qinftyalphab}
Let $p=1$ and $1\leq q\leq \infty$. Suppose that $T_{a, b, c}$ is bounded from $L^1_\alpha$ to $L_\beta^q$, then $\alpha\leq b$.
\end{lemma}

\begin{proof}
It is clear that the adjoint operator $T_{a, b, c}^\ast$, which is given by \eqref{tstar} or \eqref{tstar1}, is bounded from $(L_\beta^q)^\ast$ to $L^\infty$.
The result can be obtained by a similar argument to the  proof of Lemma~\ref{lem:aqbeta}.
\end{proof}

Let $H(\Bn)$ be the space of all holomorphic functions in $\Bn$. In the following, the notation $A^p_{\alpha}: =L^{p}_{\alpha}\cap H(\Bn)$ denotes the weighted Bergman space on $\Bn$ with the same norm as $L_\alpha^p$ for $1\leq p<\infty$.

\subsection{Proof of (ii)$\Rightarrow$(iii) for Theorem \ref{thm:pq3}}

We divide our proof into the following two cases.

\subsubsection{Case I: $1\leq q<p< \infty$ }

We need the following two lemmas to prove this case.

\begin{lemma}\cite[Theorem 30]{ZZ}\label{lem:timage}
Suppose $1\leq p<\infty$ and $b, c$ are real parameters satisfying the following conditions
\begin{itemize}
\item[(i)]$(\alpha+1)<p(b+1)$,
\item[(ii)] $c-1$ is  not a negative integer.
\end{itemize}
Then a holomorphic function $f$ in $\Bn$ belongs to $A_{\alpha+p(c-n-1-b)}^p$ if and only if
\[
f(z)=\int_{\Bn}\frac{g(w)(1-|w|^2)^b}{(1-\langle z, w\rangle)^{c}}dv(w)
\]
for some $g\in L_\alpha^p$.
\end{lemma}

\begin{lemma}\cite[Theorem 70]{ZZ}\label{lem:inculsionqp}
Let $0<q<p<\infty$. Then $A_\alpha^p\subset A_\beta^q$ if and only if
\[
\frac{1+\alpha}p<\frac{1+\beta}q,
\]
and the inclusion is strict.
\end{lemma}

We are giving the proof of (ii)$\Rightarrow$(iii) in Theorem~\ref{thm:pq3}
for the case $1\leq q<p< \infty$
in the following proposition.

\begin{proposition}
Let $1\leq q<p< \infty$. If $T_{a, b, c}$ is bounded from $L^p_\alpha$ to $L^q_\beta$,
then
\begin{align}
&-qa<\beta+1, \label{eq:parametera}\\
&(\alpha+1)<p(b+1), \label{eq:parameterb} \\
&c<n+1+a+b+\frac{1+\beta}q-\frac{1+\alpha}p.\label{eq:papameterc}
\end{align}
\end{proposition}

\begin{proof}
First, the inequality \eqref{eq:parametera} follows from Lemma \ref{lem:aqbeta},
and the inequality \eqref{eq:parameterb} follows from Lemma \ref{lem:alphap}.
Hence it is sufficient to prove \eqref{eq:papameterc}.
Assume that \eqref{eq:papameterc} is not true.
Suppose
\begin{equation}\label{eq:contrdictionc}
c\geq n+1+a+b+\frac{1+\beta}q-\frac{1+\alpha}p.
\end{equation}
From \eqref{eq:parametera} and \eqref{eq:parameterb} we get
\begin{align*}
c-1\geq n+a+b+\frac{1+\beta}q-\frac{1+\alpha}p>n-1\geq 0.
\end{align*}
It follows from Lemma~\ref{lem:timage} that
\[
T_{0, b, c}(L_\alpha^p)=A_{\alpha+p(c-n-1-b)}^p.
\]
Since \eqref{eq:contrdictionc} implies that
\[
\frac{1+\alpha+p(c-n-1-b)}p\geq a+\frac{1+\beta}q,
\]
it follows from Lemma \ref{lem:inculsionqp} that $A_{\alpha+p(c-n-1-b)}^p$ is not contained in $A_{\beta+qa}^q$.
Therefore, there exists a function $f\in L_\alpha^p$ such that $T_{0, b, c}f\notin A_{\beta+qa}^q$. Thus $T_{0, b, c}$ is unbounded from $L_\alpha^p$ to $A_{\beta+qa}^q$. Consequently, we infer that $T_{a, b, c}$ is unbounded from $L_\alpha^p$ to $L_\beta^q$. This leads to a contradiction. The proof is complete.
\end{proof}

\subsubsection{Case II: $p=\infty$ and $1\leq q<\infty$ }

For this case we need two lemmas related to generalized Lipschitz space.
For $-\infty<\beta<\infty$, the generalized Lipschitz space $\Lambda_\beta$ consists of those $f\in H(\Bn)$ for which
\[
\sup_{z\in \Bn}(1-|z|^2)^{k-\beta} |R^kf(z)|<\infty,
\]
where $k$ is any nonnegative integer with $k>\beta$ and 
\[
Rf(z)=\sum_{k=1}^n z_k\frac{\partial f}{\partial z_k}(z)
\]
is the radial derivative of $f$, see \cite[P23]{ZZ} in detail.

\begin{lemma}\cite[Theorem 17]{ZZ}\label{lem:Bloch}
Suppose $b>-1$, $c-1$ is not a negative integer and $f\in H(\Bn)$.
Then the following are equivalent.
\begin{itemize}
\item[(i)] $f\in \Lambda_{-(c-n-1-b)}$.
\item[(ii)] There exists a function $g\in L^\infty(\Bn)$ such that
\[
f(z)=\int_{\Bn}\frac{g(w)(1-|w|^2)^b}{(1-\langle z, w\rangle)^{c}}dv(w)
\]
\end{itemize}
\end{lemma}

\begin{lemma}\label{lem:growthspace}\cite[Theorem 66]{ZZ}
Let $0<p<\infty$, and let $\beta$ be any real number. Then for any $\gamma<(1+\beta)/p$ we have
\[
\Lambda_{-\gamma}\subset A_\beta^p\subset\Lambda_{-(n+1+\beta)/p}.
\]
Both inclusion are strict and best possible, where ``best possible" means that,
for each $p$ and $\beta$, the index $\gamma$ of $\Lambda_{-\gamma} $
on the left-hand side cannot be replaced by a larger number,
and the index $(n+1+\beta)/p$ on the right-hand side cannot be replaced by a smaller one.
\end{lemma}

We give the proof of (ii)$\Rightarrow$(iii) in Theorem \ref{thm:pq3}
for the case $p=\infty$ and $1\leq q<\infty$
in the following proposition.

\begin{proposition}\label{lem:inftyqn2}
Suppose $p=\infty$ and $1\leq q<\infty$. If  $T_{a, b,c}$ is bounded from  $L^\infty$ to $L_\beta^q$ then
\begin{align}
&-qa<\beta+1, \label{eq:parameterainfty}\\
&b>-1, \label{eq:parameterbinfty}\\
&c<n+1+a+b+\frac{1+\beta}q.\label{eq:parametercinfty}
\end{align}
\end{proposition}

\begin{proof}
First, \eqref{eq:parameterainfty} follows from Lemma \ref{lem:aqbeta},
and \eqref{eq:parameterbinfty} follows from Lemma \ref{lem:b1}.
Now we only need to prove \eqref{eq:parametercinfty}.
Assume that \eqref{eq:parametercinfty} is not true.
Suppose
\begin{equation}\label{eq:parcc}
c\geq n+1+a+b+\frac{1+\beta}q.
\end{equation}
Then from \eqref{eq:parameterainfty} and \eqref{eq:parameterbinfty} we get that
\[
c-1\geq n+a+b+\frac{1+\beta}q>n-1\geq0.
\]
It follows from Lemma \ref{lem:Bloch} that
\[
T_{0, b, c} (L^\infty)=\Lambda_{-(c-n-1-b)}.
\]
By \eqref{eq:parcc} we know that
\[
c-n-1-b\geq a+\frac{1+\beta}q.
\]
It follows from Lemma~\ref{lem:growthspace} that $\Lambda_{-(c-n-1-b)}$ is not contained in Bergman space $A_{\beta+aq}^q$.
Therefore, there exists a function $f\in L^\infty$ such that $T_{0, b,c}f\notin A_{\beta+aq}^q$.
Thus $T_{0, b, c}$ is unbounded from $L^\infty$ to $A_{\beta+qa}^q$.
Consequently, we get that $T_{a, b, c}$ is unbounded from $L^\infty$ to $L_\beta^q$,
which leads to a contradiction. The proof is complete.
\end{proof}

\subsection{Proof of (ii)$\Rightarrow$(iii) for Theorem \ref{thm:thm4}}
Recall that the Berezin transform  of $f\in L_\alpha^1$ is given by
\begin{align}\label{eq:berezin}
B_\alpha f(z)~&=\langle fk_z, k_z\rangle\notag\\
~&=\int_{\Bn}\frac{(1-|z|^2)^{n+1+\alpha}}
{|1-\langle z, w\rangle|^{2(n+1+\alpha)}}f(w)dv_\alpha(w), \quad z\in \Bn,
\end{align}
where
$$
k_z=K_z(w)/\sqrt{K_z(z)},  \quad w\in \Bn
$$
with
\[
K_z(w)=\frac{1}{(1-\langle w, z\rangle)^{n+1+\alpha}}
\]
is the  reproducing kernel for $A^2_{\alpha}$, see \cite[Chapter 2]{ZhuBn}. We need the following two lemmas.

\begin{lemma}\label{eq:berezint}
If $f\in H^\infty(\Bn)$, then $B_\alpha f=f$, where $H^\infty(\Bn)$ denotes the spaces of bounded holomorphic functions on $\Bn$.
\end{lemma}

\begin{proof}
Since $f\in H^\infty(\Bn)$, for each fixed $z\in \Bn$, $fk_z\in A_\alpha^2$. It follows from \cite[Theorem 2.2]{ZhuBn} that
\[
B_\alpha f(z)=\frac{1}{\sqrt{K_z(z)}}\langle fk_z, K_z\rangle=\frac{f(z)k_z(z)}{\sqrt{K_z(z)}}=f(z).
\]
\end{proof}

\begin{lemma}
Let $c\in \mathbb{R}$. Then
\begin{equation}\label{lem:berezintransform}
\int_{\Bn}\frac{(1-|\xi|^2)^{n+1+\alpha}}{(1-\langle z, w\rangle)^c|1-\langle  \xi, w\rangle|^{2(n+1+\alpha)}}dv_\alpha(w)=\frac{1}{(1-\langle z, \xi\rangle)^c}
\end{equation}
for every $z\in \Bn$.
\end{lemma}

\begin{remark}
A more general form of formula \eqref{lem:berezintransform} can be found in \cite[Lemma 2.3]{LC}.
We give a simple proof here for completion.
\end{remark}

\begin{proof}
It is clear that, for any fixed $z\in \Bn$, the function
\[
f_z(w)=\frac{1}{(1-\langle w, z\rangle)^c}\in H^\infty(\Bn).
\]
It follows from Lemma~\ref{eq:berezint} that $B_\alpha f_z=f_z$.
\end{proof}

The proof of (ii)$\Rightarrow$(iii) for Theorem \ref{thm:thm4} is given
in the following proposition.

\begin{proposition}\label{lem:tabc1infty}
Suppose $T_{a,b,c}$ is bounded from $L_\alpha^1$ to $L^\infty$, then
\[
\bigg\{
\begin{aligned}
&a\geq 0,\  \alpha\leq b, \\
&c\leq a+b-\alpha.
\end{aligned}
\]
\end{proposition}

\begin{proof}
It follows from Lemma~\ref{lem:pinfty1a} that  $a\geq 0$,
and it follows from Lemma~\ref{lem:qinftyalphab} that  $\alpha\leq b$.
Thus we only need to prove that $c\leq a+b-\alpha$.
We consider two cases.

\emph{Case I: $\alpha=b$.}
It is sufficient for us to show that $c\leq a$ in this case.
For any fixed $\xi\in \Bn$, consider the function
\[
f_{\xi}(z)=\frac{(1-|\xi|^2)^{n+1+\alpha}}{|1-\langle z, \xi\rangle|^{2(n+1+\alpha)}}.
\]
It follows from Lemma~\ref{lem:formula} that there is a constant $C$,
independent of $\xi$, such that $\|f_\xi\|_{1, \alpha}\leq C$.
Applying \eqref{lem:berezintransform} we get
\begin{align*}
 T_{a, b, c}f_\xi(z)~&=(1-|z|^2)^a\int_{\Bn}\frac{(1-|w|^2)^\alpha}{(1-\langle z, w\rangle)^c}f_\xi(w)dv(w)\\
~&=\frac{1}{c_{\alpha}}\,\frac{(1-|z|^2)^a}{(1-\langle z, \xi\rangle)^c}.
\end{align*}
Since $T_{a,b,c}$ is bounded from $L_\alpha^1$ to $L^\infty$, it follows that
\[
\|T_{a, b, c}f_\xi\|_{\infty}\leq \|T_{a, b, c}\|_{L_\alpha^1\to L^\infty}\|f_\xi\|_{1, \alpha}
\]
for any fixed $\xi\in \Bn$. Thus we get
\[
\sup_{\xi\in \Bn}\bigg\|\frac{(1-|z|^2)^a}{(1-\langle z, \xi\rangle)^c}\bigg\|_{\infty}<+\infty,
\]
which implies that $c\leq a$.

\emph{Case II: $\alpha<b$.} For any fixed $\xi\in \Bn$, consider the function
\[
f_{\xi}(z)=\frac{(1-|\xi|^2)^{b-\alpha}}{(1-\langle z, \xi\rangle)^{n+1+b}}
\]
It follows from Lemma \ref{lem:formula} that there is a constant $C$,
independently of $\xi$, such that $\|f_{\xi}\|_{1,\alpha}\leq C$.
A similar discussion as in the proof of \cite[Lemma 9]{Zh} shows that
\[
T_{a, b,c}f_\xi(z)
=\frac{1}{c_{\alpha}}\frac{(1-|z|^2)^a(1-|\xi|^2)^{b-\alpha}}{(1-\langle z, \xi\rangle)^c}.
\]
Since $T_{a,b,c}$ is bounded from $L_\alpha^1$ to $L^\infty$, we get that
\[
\sup_{\xi\in \Bn}\bigg\|\frac{(1-|z|^2)^a(1-|\xi|^2)^{b-\alpha}}
{(1-\langle z, \xi\rangle)^c}\bigg\|_\infty<+\infty,
\]
which implies that $c\leq a+b-\alpha$. The proof is complete.
\end{proof}

\subsection{Proof of (ii)$\Rightarrow$(iii) for Theorem \ref{thm:thm5}}

The following two propositions are devoted to giving the proof of (ii)$\Rightarrow$(iii) for Theorem \ref{thm:thm5}.

\begin{proposition}\label{lem:necessitypinfty}
Let $1< p<\infty$ and $q=\infty$.  If  $T_{a, b,c}$ is bounded from $L_\alpha^p$ to $L^\infty$ then
\begin{equation}\label{eq:abc2}
\bigg\{
\begin{aligned}
&a=0,  \quad \alpha+1<p(b+1)\\
&c<n+1+b-\frac{n+1+\alpha}p.
\end{aligned}
\end{equation}
or
\begin{equation}\label{eq:abc1}
\bigg\{
\begin{aligned}
&a>0, \quad (\alpha+1)<p(b+1) \\
&c\leq n+1+a+b-\frac{n+1+\alpha}p,
\end{aligned}
\end{equation}
\end{proposition}

\begin{proof}
Since $T_{a, b,c}$ is bounded from $L_\alpha^p$ to $L^\infty$,
and $L^1$ injects continuously into $(L^\infty)^*$,
it follows that its adjoint operator $T_{a, b,c}^\ast$ given in \eqref{tstar1} is bounded from $L^1$ to $L_\alpha^{p^\prime}$.
Applying Theorem~\ref{thm:pq2} we get
\[
\bigg\{
\begin{aligned}
&a=0,\quad -p^\prime(b-\alpha)<\alpha+1,\\
&c<b-\alpha+\frac{n+1+\alpha}{p^\prime},
\end{aligned}
\]
and
\[
\bigg\{
\begin{aligned}
&a>0,\quad -p^\prime (b-\alpha)<\alpha+1,\\
&c\leq b-\alpha+a+\frac{n+1+\alpha}{p^\prime},
\end{aligned}
\]
which are equivalent to \eqref{eq:abc2} and \eqref{eq:abc1}. The proof is complete.
\end{proof}

\begin{proposition}\label{lem:necessityinftyinfty}
Suppose $T_{a,b,c}$ is bounded from $L^\infty$ to $L^\infty$, then
\[
\bigg\{
\begin{aligned}
&a=0, \ b>-1, \\
&c<n+1+b.
\end{aligned}
\quad
\mathrm{or}
\quad
\bigg\{
\begin{aligned}
&a>0, \ b>-1, \\
&c\leq n+1+a+b.
\end{aligned}
\]
\end{proposition}

\begin{proof}
First, it follows from Lemma \ref{lem:b1} that $b>-1$, and it follows from Lemma \ref{lem:pinfty1a} that $a\geq 0$.
Assume that $ c\geq n+1+b$ when $a=0$, and  $c>n+1+a+b$ when $a>0$.
Recall that, in this case the adjoint operator of $T_{a,b,c}$ is given by
\[
T_{a,b,c}^{\ast}f(z)=(1-|z|^2)^b\int_{\Bn}\frac{(1-|w|^2)^a}{(1-\langle z, w\rangle)^c}f(w)dv(w).
\]
Applying Theorem B we have that $T_{a,b,c}^\ast$ is unbounded from $L^1$ to $L^1$.
This implies that $T_{a,b,c}$ is unbounded from $L^\infty$ to $L^\infty$, which is a contradiction.
The proof is complete.
\end{proof}

\section{Applications}

In this section, we first apply our results to the following operator,
\[
K_c^\alpha f(z)=\int_{\Bn}\frac{f(w)}{(1-\langle z, w\rangle)^c}\,dv_{\a}(w),
\]
which was first considered in \cite{ChFangWY} for the case of the unit disk with $\a=0$,
and generalize  Theorems 1-3 in \cite{ChFangWY} to weighted Lebesgue spaces on the unit ball of $\C^n$.

The following result, for the case of the unit disk with $\a=0$,
was first proved in \cite[Theorem 1]{ChFangWY}.

\begin{corollary}\label{cor:cc}
Let $-1<\alpha<\infty$.
Then there exists a pair $(p, q)\in [1,+\infty]\times[1,+\infty]$
such that $K_c^\alpha: L_\alpha^p\to L_\alpha^q$ is bounded if and only if $c<n+2(1+\alpha)$.
\end{corollary}

\begin{proof}
We first prove the ``only if" part.
Assume that there exists a pair $(p, q)\in [1, +\infty]\times[1, +\infty]$
such that $K_c^\alpha: L_\alpha^p\to L_\alpha^q$ is bounded.
Examining the inequalities for $c$ in Theorems~A, B
and Theorems \ref{thm:pq3}-\ref{thm:thm5} with $a=0$, $b=\a$ and $\beta=\a$,
we can see that, for each case we have
\[
c<n+1+\alpha+\frac{1+\alpha}q.
\]
If $q=1$ then the right-hand side reaches the maximum $n+2(1+\alpha)$.
Thus we must have $c<n+2(1+\alpha)$.

Next we prove the ``if" part.
Assume that $c<n+2(1+\alpha)$.
The proof of this part will be divided into three cases.

\emph{Case I: $c\leq 0$.}
It is obviously that $K_c^\alpha: L_\alpha^p\to L_\alpha^q$ is bounded
for all $(p, q)\in [1, +\infty]\times[1, +\infty]$ for this case
by Theorems~A, B and Theorems~\ref{thm:pq3}-\ref{thm:thm5}.

\emph{Case II: $0<c\leq n+1+\alpha$.}
Let $\tau_1=n+1+\alpha-c\geq 0$. Then we can choose $p, q$ such that $1<p\leq q<\infty$ and
\[
\frac{1}{p}-\frac{1}q=\frac{\tau_1}{n+1+\alpha}.
\]
Thus
\[
c=n+1+\alpha-\tau_1=n+1+\alpha+\frac{n+1+\alpha}q-\frac{n+1+\alpha}p.
\]
By Theorem A, we get that $K_c^\alpha: L_\alpha^p\to L_\alpha^q$ is bounded.

\emph{Case III: $n+1+\alpha<c<n+2(1+\alpha)$.}
Let $\tau_2=n+2(1+\alpha)-c$.
Since $n+1+\a<c$ we have $0<\tau_2<1+\alpha$.
Thus we can choose $p, q$ such that $1\leq q<p<\infty$ and
\[
\frac 1q-\frac 1p>1-\frac{\tau_2}{1+\alpha}.
\]
This implies that
\begin{align*}
c=n+2(1+\alpha)-\tau_2<n+1+\alpha+\frac{1+\alpha}q-\frac{1+\alpha}p.
\end{align*}
It follows from Theorem~\ref{thm:pq3} that $K_c^\alpha: L_\alpha^p\to L_\alpha^q$ is bounded.
The proof is complete.
\end{proof}

The following result was first obtained as \cite[Theorem 2]{ChFangWY}
for the case of the unit disk and $\a=0$.

\begin{corollary}\label{cor:corleq}
Let $0<c\leq n+1+\alpha$ and  $1\leq p, q\leq \infty$.
Then $K_c^\alpha: L_\alpha^p\to L_\alpha^q$ is bounded if and only if $p, q$ satisfy one of the following conditions:
\begin{itemize}
\item[(i)] $p=1,\  q<\dfrac{n+1+\alpha}c. $
\item[(ii)] $1<p<\dfrac{n+1+\alpha}{n+1+\alpha-c}$, and
\begin{equation}\label{eq: pqgeq}
 \frac 1q\geq \frac 1p+\frac c{n+1+\alpha}-1.
\end{equation}
\item[(iii)]$p=\dfrac{n+1+\alpha}{n+1+\alpha-c},\ q<\infty$.
\item[(iv)] $p>\dfrac{n+1+\alpha}{n+1+\alpha-c}$.
\end{itemize}
\end{corollary}

\begin{proof}
We divide the proof into four steps.

\emph{Step I: Assume that $p=1$.} It is clear that $q<(n+1+\alpha)/c$ is equivalent to
\[
c<\frac{n+1+\alpha}q.
\]
Then, by Theorem B, we get that $K_c^\alpha: L_\alpha^p\to L_\alpha^q$
is bounded if and only if  $p, q$ satisfy condition (i).

\emph{Step II: Assume that $1<p<(n+1+\alpha)/(n+1+\alpha-c)$.}
Then we get
\begin{equation}\label{eq:c1alpha}
c>n+1+\alpha-\frac{n+1+\alpha}p.
\end{equation}
It follows from Theorem \ref{thm:thm5} that $K_c^\alpha: L_\alpha^p\to L_\alpha^q$ is unbounded for $q=\infty$.
On the other hand, from 
\eqref{eq: pqgeq} and \eqref{eq:c1alpha} we get
\[
\frac1q\geq\frac 1p+\frac c{n+1+\alpha}-1>0,
\]
which implies that $q\neq \infty$.
Hence we only need to consider the case $q<\infty$.
It suffices for us to consider the following two cases.

\emph{Case I: $1<p\leq q<\infty$.} By theorem A, we obtain that $K_c^\alpha: L_\alpha^p\to L_\alpha^q$ is bounded if and only if
\begin{equation}\label{eq:cnalpha}
c\leq n+1+\alpha+\frac{n+1+\alpha}q-\frac{n+1+\alpha}p,
\end{equation}
which is equivalent to \eqref{eq: pqgeq}.

\emph{Case II:  $1\leq q<p<(n+1+\alpha)/(n+1+\alpha-c)$.}
In this case we have
$$
\frac1q>\frac1p\ge\frac 1p+\frac c{n+1+\alpha}-1,
$$
and so \eqref{eq: pqgeq} always holds.
Thus we only need to show that $K_c^\alpha: L_\alpha^p\to L_\alpha^q$ is bounded for all $p, q$ in this case.
Since $c\leq n+1+\alpha$ and $q<p$, we get that
\begin{equation}\label{eq:ccondtion2}
c<n+1+\alpha+\frac{1+\alpha}q-\frac{1+\alpha}p.
\end{equation}
By Theorem \ref{thm:pq3}, we see that $K_c^\alpha: L_\alpha^p\to L_\alpha^q$ is bounded.

\emph{Step III. Assume that $p=(n+1+\alpha)/(n+1+\alpha-c)$.}  Then we have
\[
c=n+1+\alpha-\frac{n+1+\alpha}p.
\]
It follows from Theorem \ref{thm:thm5} that $K_c^\alpha: L_\alpha^p\to L_\alpha^q$ is unbounded for $q=\infty$,
and it also follows from Theorem A and Theorem \ref{thm:pq3} that  $K_c^\alpha: L_\alpha^p\to L_\alpha^q$ is bounded for $q<\infty$
since $c$ satisfies both \eqref{eq:cnalpha} and \eqref{eq:ccondtion2}.

\emph{Step IV.  Assume that $p>(n+1+\alpha)/(n+1+\alpha-c)$.} Then we get
\[
c<n+1+\alpha-\frac{n+1+\alpha}p.
\]
It follows Theorem A,  Theorem~\ref{thm:pq3} and Theorem~\ref{thm:thm5}
that $K_c^\alpha: L_\alpha^p\to L_\alpha^q$ is bounded for all $1\leq q\leq \infty$.
The proof is complete.
\end{proof}

The following result was first obtained as \cite[Theorem 3]{ChFangWY}
for the case of the unit disk and $\a=0$.

\begin{corollary}\label{cor:ccgeq}
Let $n+1+\alpha<c< n+2(1+\alpha)$ and  $1\leq p, q\leq \infty$. Then $K_c^\alpha: L_\alpha^p\to L_\alpha^q$ is bounded if and only if $p, q$ satisfy
\[
\frac{1+\alpha}{n+2(1+\alpha)-c}<p\leq \infty, \quad \frac 1q>\frac 1p+\frac{c-(n+1+\alpha)}{1+\alpha}.
\]
When $p=\infty$, the second inequality should be interpreted as
\[
\frac 1q>\frac{c-(n+1+\alpha)}{1+\alpha}.
\]
\end{corollary}
\begin{proof}
Since $n+1+\alpha<c< n+2(1+\alpha)$, it follows from Theorem A, Theorem B and Theorem~\ref{thm:thm5}
that $K_c^\alpha: L_\alpha^p\to L_\alpha^q$ is unbounded if $p, q$ satisfy $1\leq p\leq q\leq \infty$.
Thus we only need to discuss the case when $p, q$ satisfy $1\leq q<p\leq \infty$.
It follows from Theorem~\ref{thm:pq3} that $K_c^\alpha: L_\alpha^p\to L_\alpha^p$ is bounded
if and only if $c$ satisfies \eqref{eq:ccondtion2}, which is equivalent to
\[
\frac 1q>\frac 1p+\frac{c-(n+1+\alpha)}{1+\alpha}.
\]
Since $1/q\leq 1$, from the above inequality we must have
\[
\frac 1p+\frac{c-(n+1+\alpha)}{1+\alpha}<1,
\]
which is equivalent to
\[
p>\frac{1+\alpha}{n+2(1+\alpha)-c}.
\]
The proof is complete.
\end{proof}

Finally we specify another two important special cases. One is the weighted Bergman projection.
Recall that for $\gamma>-1$, the weighted Bergman projection is given by the following formula.
\[
P_\gamma f(z)=\int_{\Bn}\frac{f(w)}{(1-\langle z, w\rangle)^{n+1+\gamma}}dv_{\gamma}(w).
\]
The other one is the Berezin transform $B_\gamma$ defined in \eqref{eq:berezin}.
It is clear that $P_{\gamma}=c_{\gamma}T_{0,\gamma,n+1+\gamma}$ and
$B_{\gamma}=c_{\gamma}S_{n+1+\gamma,\gamma,2(n+1+\gamma)}$.
Hence, as special cases of our main results,
we have the following characterizations of
the $L^p$-$L^q$ boundedness of $P_\gamma$ and $B_\gamma$.

\begin{corollary}\label{cor:pr}
Suppose $1\leq p, q\leq \infty$.
\begin{itemize}
\item[(i)] When $1<p\leq q<\infty$, $P_\gamma$ (or $B_\gamma$) is bounded from $L_\alpha^p$ to $L_\beta^q$ if and only if
\[
\left\{
\begin{aligned}
&\alpha+1<p(\gamma+1)\\
&(n+1+\alpha)/p\leq (n+1+\beta)/q.
\end{aligned}
\right.
\]
\item[(ii)]
When $1=p\leq q<\infty$, $P_\gamma$ (or $B_\gamma$) is bounded from $L_\alpha^1$ to $L_\beta^q$ if and only if
\[
\left\{
\begin{aligned}
&\alpha<\gamma\\
& n+1+\alpha\leq (n+1+\beta)/q
\end{aligned}
\right.
\quad
\textrm{or}
\quad
\left\{
\begin{aligned}
& \alpha=\gamma\\
& n+1+\alpha<(n+1+\beta)/q
\end{aligned}
\right.
\]
\item[(iii)]When $1\leq q<p<\infty$, $P_\gamma$ (or $B_\gamma$) is bounded from $L_\alpha^p$ to $L_\beta^q$ if and only if
\[
\left\{
\begin{aligned}
&(\alpha+1)<p(\gamma+1)\\
&(1+\alpha)/p<(1+\beta)/q
\end{aligned}
\right.
\]
\item[(iv)] $P_\gamma$ is a bounded from $L^\infty$ to $L_\beta^q$ if and only if $1\leq q<\infty$,
and $B_\gamma$ is bounded from $L^\infty$ to $L_\beta^q$ if and only if $1\leq q\leq \infty$.
\end{itemize}
\end{corollary}

\section*{Acknowledgements}
The authors would like to thank Guangyao Zhang for his help of drawing figures.

\section*{Additional Remarks after publication}

This paper was originally published in [Journal of Functional Analysis, 282 (2022), 109345, 
DOI: https://doi.org/10.1016/j.jfa.2021.109345].
After this paper was published in 2021, the authors found that 
Kaptanogl\v{u} and \"{U}reyen \cite{KU} obtained some similar results 
earlier for integral operators with Bergman-Besov kernels, which when 
$\alpha>-1$ overlap some of our main results here. However, we note that the proofs in this 
paper and that of \cite{KU} are different. For sufficiency, we used different 
versions of Schur’s tests (for example, see Lemma~\ref{lem:1} 
in this paper and Theorem 5.2 in \cite{KU}); 
and more essentially, for necessity, our proofs here employ the 
integral representations in Lemma~\ref{lem:timage} and Lemma~\ref{lem:Bloch}, 
while in \cite{KU} it uses a versatile method depending on precise 
imbeddings of certain spaces of analytic functions.

\end{document}